\documentclass[a4paper, 10 pt, conference]{ieeeconf}%
\usepackage{graphics}
\usepackage{epsfig}
\usepackage{mathptmx}
\usepackage{amsmath}
\usepackage{amssymb}
\usepackage{setspace}
\usepackage[colorinlistoftodos]{todonotes}
\usepackage[colorlinks=true, allcolors=blue]{hyperref}
\usepackage{caption}
\usepackage{subcaption}
\usepackage{hyperref}
\usepackage{xcolor,colortbl}
\usepackage{setspace}
\usepackage{float}
\usepackage{cite}
\usepackage{graphicx}
\usepackage[space]{grffile}
\usepackage{amsfonts}%
\usepackage{amsmath}
\usepackage{booktabs}
\usepackage{xcolor,colortbl}
\usepackage{optidef}
\usepackage{tikz}
\usepackage{breqn}
\usepackage[font={small}]{caption}
\PassOptionsToPackage{usenames,dvipsnames,svgnames}{xcolor}
\usetikzlibrary{arrows,positioning,automata}
\usepackage{smartdiagram}
\usetikzlibrary{arrows,positioning,automata}
\providecommand{\U}[1]{\protect\rule{.1in}{.1in}}
\IEEEoverridecommandlockouts
\begin{document}
\bstctlcite{MyBSTcontrol}

\title{{\LARGE \textbf{The Penetration Rate Effect of Connected and Automated Vehicles in Mixed Traffic Routing }}}
\author{Arian Houshmand$^{1}$, Salom\'{o}n Wollenstein-Betech$^{1}$, and Christos G. Cassandras$^{1}$\thanks{*This work was 
supported in part by NSF under grants ECCS-1509084, DMS-1664644,
CNS-1645681, by AFOSR under grant FA9550-15-1-0471, by ARPA-E’s
NEXTCAR program under grant DE-AR0000796 and by the MathWorks.}\thanks{$^{1}$The authors are with the
Division of Systems Engineering, Boston University, Brookline, MA 02446 USA
\texttt{{\small arianh@bu.edu; salomonw@bu.edu; cgc@bu.edu}}}}
\maketitle
\begin{abstract}
We study the problem of routing  Connected and Automated Vehicles (CAVs) in the presence of mixed traffic (coexistence of regular vehicles and CAVs). In  this setting, we assume that all CAVs belong to the same fleet, and can be routed using a centralized controller. The routing objective is to  minimize a given overall fleet traveling cost (travel time or energy consumption). We assume that  regular  vehicles (non-CAVs)  choose  their  routing decisions  selfishly  to  minimize their  traveling  time.  We  propose  an algorithm  that  deals  with  the routing interaction  between  CAVs and  regular uncontrolled vehicles.  We  investigate  the effect of assigning system-centric  routes under  different  penetration  rates  (fractions) of CAVs. To validate our method, we apply the proposed routing algorithms to the Braess Network and to a sub-network of the Eastern Massachusetts (EMA) transportation network using actual traffic data provided by the Boston Region Metropolitan Planning Organization. The results suggest that collaborative routing decisions of CAVs improve not only the cost of CAVs, but also that of the non-CAVs. Furthermore, even a small CAV penetration rate can ease congestion for the entire network.

\end{abstract}

\section{INTRODUCTION}
Every year Americans face more than 6.9 billion hours of delay in traffic which costs the US more than 160 billion dollars in urban congestion costs \cite{david_schrank_bill_eisele_tim_lomax_jim_bak_2015_2015}. In addition, due to heavy traffic congestion, an annual amount of 3.1 billion gallons of fuel is being wasted in traffic \cite{david_schrank_bill_eisele_tim_lomax_jim_bak_2015_2015}.

The advent of Connected and Automated Vehicles (CAVs) has been facilitated by the emergence of vehicle automation technologies, as well as new forms of telecommunication technologies, such as Dedicated Short-Range Communication (DSRC) \cite{kenney_dedicated_2011} and 5G \cite{andrews_what_2014}. The latter has enabled Vehicle-to-Vehicle (V2V) and Vehicle-to-Infrastructure (V2I) communication capabilities. Therefore, CAVs can help reduce traffic congestion and environmental impacts of our daily commute, as well as improve safety through  collaborative decisions.   

\begin{figure}[ptb]
\centering
 \includegraphics[width=0.80\linewidth]{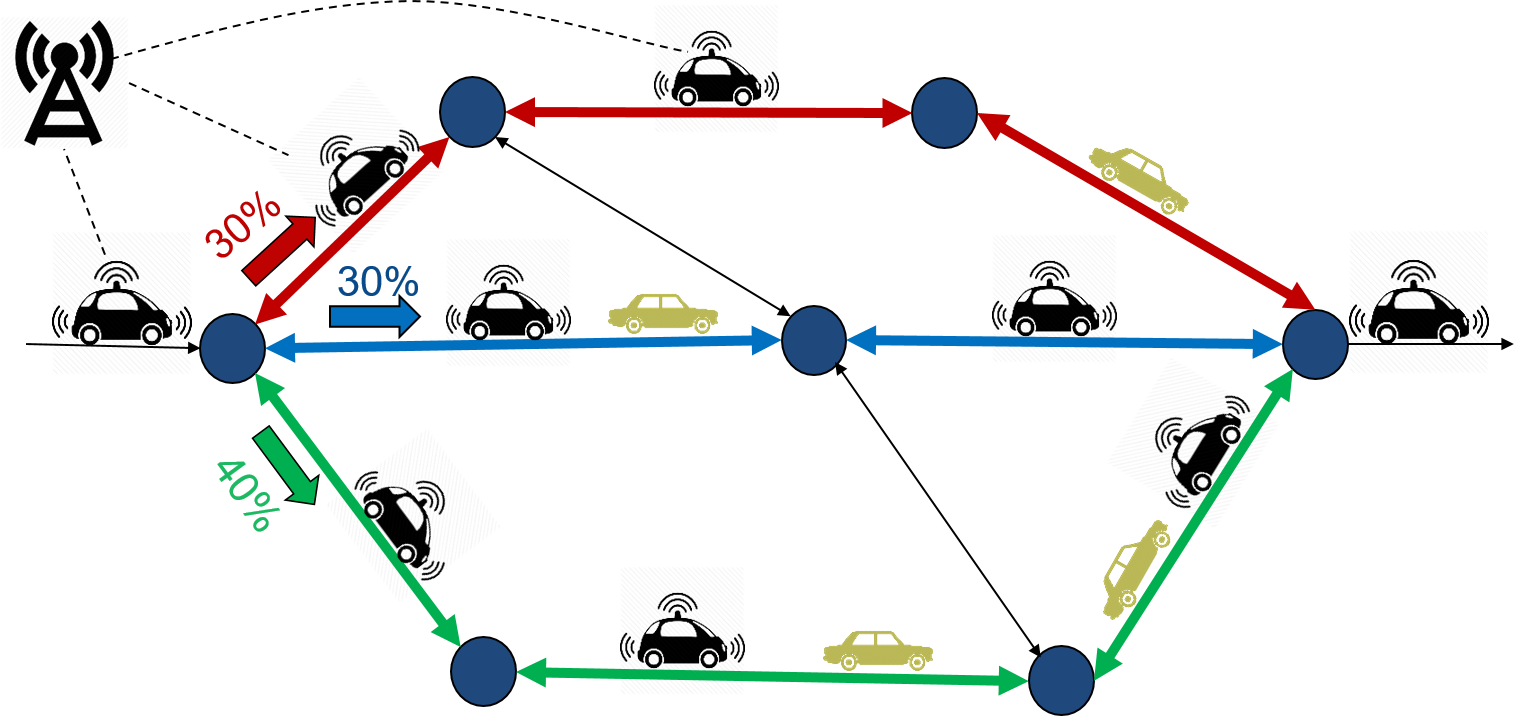}
  \caption{The centralized controller is assigning routes to CAVs entering  the network. Red, blue, and green links show three different routes for CAVs.}  
  \label{fig: CAV schematic}
\end{figure}

Many studies have been performed to investigate how CAVs can transform the future of cities \cite{guanetti_control_2018}. For instance, we may be able to eliminate traffic lights and create unsignalized intersections to reduce congestion and energy consumption \cite{malikopoulos_decentralized_2018}. Another interesting area of focus is cooperative adaptive cruise control (CACC)  \cite{van_arem_impact_2006,rajamani_experimental_2001}, an extension of adaptive cruise control (ACC). By  exploiting V2V communication, CACC can reduce the headway between vehicles. As a result, CACC can reduce drag forces and increase energy efficiency \cite{al_alam_experimental_2010}, as well as road capacity and throughput \cite{lazar_capacity_2017,van_arem_impact_2006}. 


Additionally, an increasing number of autonomous fleet operated businesses are emerging, namely autonomous mobility on demand (AMoD) \cite{gibbs_google_2017}, delivery services, and freight shipping.
Coming up with strategies to minimize the travel time and energy cost of these fleets not only reduces traffic congestion and carbon emission impacts, but also increases profitability.

In this paper, we seek to find how optimizing routing decisions of CAVs affects the overall energy consumption costs and travel times of all vehicles. We investigate the interaction between CAVs and regular vehicles and their effects on total travel time and energy consumption in traffic networks. We assume that all CAVs belong to the same fleet (e.g., AMoD), and the fleet operator is trying to minimize their costs (energy or time) by systematically routing the fleet given their origin-destination (O-D) demand. 

Similar to this work, Mehr et al. \cite{mehr_can_nodate} studied how the presence of CAVs can affect  mobility in traffic networks. They assumed  CAVs can benefit from CACC by creating shorter headway which increases the road capacity. In this context, they adopted the mixed traffic road capacity model   from \cite{lazar_capacity_2017}. Their results show that if all vehicles (CAVs and non-CAVs) make selfish routing decisions, the presence of CAVs might worsen traffic conditions. 
In contrast to \cite{mehr_can_nodate}, we adopt the viewpoint that there is a centralized controller capable of routing all CAVs given their origin and destination. Moreover, we do not assume the shorter headway for CAVs which was considered in \cite{mehr_can_nodate}. We show that optimal routing of CAVs under these assumptions can not only benefit CAVs, but also help non-CAVs to save time and energy.

The contributions of this paper are summarized as follows. We first review a system-centric (socially optimal) routing algorithm that minimizes the total travel time assuming 100\%  CAVs in the system. We then propose  algorithms which assign system-centric time-optimal or energy-optimal routes (eco-routes) to CAVs in the presence of mixed traffic (both CAVs and non-CAVs in the system). 
The eco-routing algorithm can handle  different vehicle classes including electric vehicles (EVs), hybrid electric vehicles (HEVs), plug-in hybrid electric vehicles (PHEVs), and conventional vehicles (CVs). Additionally, using the notion of  \emph{Wardrop equilibrium} \cite{beckmann_studies_1955}, we model the user-centric (selfish) routing decisions of non-CAVs. 

Our results indicate that optimal routing of CAVs can benefit both CAVs and non-CAVs in energy savings and travel times. In addition, we study the performance of the routing algorithms for various CAV penetration rates. We provide evidence that even under small CAV penetration rates, CAVs and non-CAVs benefit. 

The remainder of this paper is organized as follows. In Section \ref{sec: time optimal routing} we propose an algorithm which assigns system-centric time-optimal routes to CAVs in the presence of mixed traffic. In Section \ref{sec: mixed traffic energy}, after reviewing different energy models for various vehicle classes, we formulate the system-centric eco-routing problem for CAVs to minimize their overall energy consumption costs in mixed traffic. In Section \ref{sec: non-CAV flow modeling}, we review the Traffic Assignment Problem and propose a framework to model non-CAV flow. In Section \ref{sec: numerical results}, we use both a simple example and actual historical data to validate the performance of our routing algorithms. Finally, conclusions and further research directions are outlined in Section \ref{sec: conclusions}.


\section{Time-Optimal Routing} \label{sec: time optimal routing}
The objective of the Time-Optimal routing  problem is to minimize the overall travel time of CAVs. To achieve this goal, we assume (1) the central controller for CAVs has full information on the Origin-Destination (O-D) demand of both CAVs and non-CAVs and (2)  non-CAVs route themselves selfishly (i.e., use the route that minimizes their individual travel time).
In order to build a model that captures the effect of mixed traffic on the optimal routing of CAVs, we first consider all vehicles to be CAVs and that they can be controlled centrally. In this respect, in Section \ref{sec: system centric time optimal}, we calculate the system-centric (social) optimal solution for the 100\% CAV penetration rate. Subsequently, in Section \ref{sec: mixed traffic time},  we generalize the routing model to find optimal routes for CAVs in mixed traffic scenarios.

\subsection{System-Centric Time-Optimal Routing }\label{sec: system centric time optimal}
 First we assume an all-CAV network, and we can route them using a centralized controller. The system-centric objective is to minimize  total traveling time (delay) of CAVs in the network.  In particular, we seek to find the route occupancy matrix (probabilities) for allocating vehicles to routes. In other words, find the probability matrix $\textbf{P}$ where its elements denote the probability that a vehicle traveling from an origin \textit{O} to  destination \textit{D} uses route $r$. 
\subsubsection{Problem Formulation} \label{sec: problem formulation 100 CAV time} 
 As in \cite{zhang_price_2018}, we model the traffic network as a directed graph $G=(\mathcal{V},\mathcal{A},\mathcal{W})$
where $\mathcal{V}$ is the set of nodes, $\mathcal{A}$ is the set of links, and \(\mathcal{W}=\{\textbf{w}_{i}:\textbf{w}_{i}=(w_{si},w_{ti}),i \in [\![\mathcal{W}]\!]\}\) is the set of all O-D pairs. We assume that all O-D pairs start and end at one of the network's nodes. Let the node-link incidence matrix for the strongly-connected and directed graph  $G$ be denoted by \(\textbf{N}\in \{0,1,-1\}^{|\mathcal{V}|\times|\mathcal{A}|}\), and let the link-route incidence matrix be denoted by \textbf{A}. Let us define \(d^{w}\geq0\) as the flow demand from \(w_{s}\) to \(w_{t}\) for any O-D pair \(\textbf{w}=(w_{s},w_{t})\). Moreover, the route choice probability matrix is defined as \(\textbf{P}=[p_{ir}]\), where \(p_{ir}\) is the probability of taking route \textit{r} while traveling through O-D pair \textit{i}. Let \(\textbf{g}=(g_{i}; i\in [\![\mathcal{W}]\!])\) be the O-D demand vector. 

Let us define the power-set of routes  \(\mathcal{\textbf{R}}=\{\mathcal{R}_{i} ;i\in [\![\mathcal{W}]\!]\}\), where \(\mathcal{R}_{i}\) is the set of allowable routes for each O-D pair \textit{i}. Finally, the link-route incidence matrix is denoted by \(\textbf{A}=\{\alpha_{a,r}^{i}; i\in [\![\mathcal{W}]\!], r\in \mathcal{R}_{i}, a\in \mathcal{A} \}\) in which:
\[\alpha_{a,r}^{i}=\begin{cases}
1; & \text{if route \(r\in \mathcal{R}_{i}\) uses link a}\\
0; & \text{otherwise.}
\end{cases}\]
Additionally, let \(A_{i}\) be the sub-matrix of \(\textbf{A}\) which includes the columns of \(\textbf{A}\) where \(r\in \mathcal{R}_{i}\).
The total flow is denoted by $\textbf{x}=\{x_{a};a\in \mathcal{A}\}$ where $x_{a}$ is the flow on each link $a\in \mathcal{A}$.

Considering \(a\in \mathcal{A}, i\in [\![\mathcal{W}]\!], r\in \mathcal{R}_{i}\) we can formulate the system-centric time-optimal problem as follows:
\begin{subequations} \label{eqn: 100 CAV main problem}
\begin{gather}
\label{eqn: system-centric time cost}
\min_{\textbf{P}}\sum_{a\in \mathcal{A}}t_{a}(x_{a})x_{a}
\\
\label{eqn: total flow on each link}
\textbf{x}=\textbf{AP}^{T}\textbf{g}
\\
\label{eqn: travel time func}
t_{a}(x_{a})=t_{a}^{0}\sum _{i=1}^{n}\beta_{i}(\frac{x_{a}}{m_{a}})^{(i-1)}
\\
\label{eqn: probability const. Social time}
\begin{array}
[c]{lr}
\sum_{r\in \mathcal{R}_{i}}p_{ir}=1;& \forall i\in [\![\mathcal{W}]\!]
\end{array}
\\
\begin{array}
[c]{lr}
 p_{ir}\in [0,1];&\forall i\in [\![\mathcal{W}]\!],\forall r\in \mathcal{R}_{i}
\end{array}
\end{gather}
\end{subequations}
where $t_{a}(x_{a})$ is the traveling time of link $a$ as a function of its corresponding traffic flow $x_{a}$, which can be modeled as an increasing polynomial function using  \eqref{eqn: travel time func}. $t_{a}^{0}$, and $m_{a}$  are the free flow travel time and flow capacity of link \(a \in \mathcal{A}\) respectively. Moreover, \(\boldsymbol{\beta}=(\beta_{i}, i=1,2,...,n) \) is the vector of coefficient factors for calculating traveling time in  \eqref{eqn: travel time func}. A common value is \(\boldsymbol{\beta}=\{1,0,0,0,0.15\} \) which is the US Bureau of Public Roads (BPR) travel time function \cite{manual_bureau_1964, huang_combined_1995}. The constraint \eqref{eqn: probability const. Social time} enforces the requirement that the sum of all the fractions of vehicles traveling through an O-D pair is 1.  The decision variable is the routing-probability matrix \(\textbf{P}=[p_{ir}]\) that for each O-D pair \( i\in [\![\mathcal{W}_i]\!]\) assigns fractions of vehicles to allowable existing routes \(r\in \mathcal{R}_{i}\)  between any given O-D pair. The inputs to the problem are the link-route incidence matrix ($\textbf{A}$), and the O-D demand vector $\textbf{g}$.  Note that instead of solving for individual links to follow for each vehicle, we are here assigning CAVs routes to follow between each OD pair. This transformation helps us reduce the decision space to select routes between OD pairs rather than finding link-based decisions.
The solution of Problem \ref{eqn: 100 CAV main problem} is often referred to as the system-centric or social optimal solution in the transportation literature.

\subsection{System-Centric Time-Optimal Routing in the Presence of Mixed Traffic}\label{sec: mixed traffic time}
In this section, we address the system-centric time-optimal routing of CAVs in the presence of mixed traffic (CAVs and non-CAVs). In this case, only a portion of vehicles are CAVs and can be controlled through a centralized controller. As a result, instead of finding a routing scheme that minimizes total costs for all vehicles in the system, we focus on minimizing travel time for the CAV portion of traffic. We consider all CAVs as belonging to the same fleet (e.g., AMoD) and that the fleet management company is trying to minimize total traveling time of the fleet. In order to solve this problem we make four assumptions: (1) The non-CAV traffic flow equilibrium is inferred from data (more details in Sec. \ref{sec: non-CAV flow modeling}). (2) There exists a centralized controller which can route the CAV portion of traffic. (3) Up to $m$ number of routes are chosen for every O-D pair. (4) Travel time functions $t(\cdot)$ are strongly monotone and continuously differentiable.

 \subsubsection{Problem Formulation} \label{sec: mixed CAV time}
 Let us define the CAV penetration rate $\gamma$, as the portion of traffic that consists of CAVs (fraction of CAVs in the system). As in the system-centric case for 100\% CAV, we define \(\textbf{P}_{c}=[p_{ir}^{c}]\) to be the route choice probability matrix for the CAV portion of traffic. Moreover, \(\textbf{g}_{c}=\{g_{i}^{c}; i\in [\![\mathcal{W}]\!]\}\) is the O-D demand vector for CAVs. As mentioned before, we assume that the non-CAV traffic flow equilibrium is inferred from data, and is known. Let us define $\textbf{x}^{c}=\{x_{a}^{c};a\in \mathcal{A}\}$ and $\textbf{x}^{nc}=\{x_{a}^{nc};a\in \mathcal{A}\}$ as the flow of CAVs and non-CAVs in the system respectively, where $x_{a}^{c}$ and $x_{a}^{nc}$ are the CAV and non-CAV flow on each link $a\in \mathcal{A}$.
 As a result, using the same notation as in Section \ref{sec: problem formulation 100 CAV time}, the optimization problem can be written as:
 \begin{subequations}\label{eqn: mixed time main problem}
 \begin{gather}
\label{eqn: mixed CAV time  objective}
\min_{\textbf{P}_c}\sum_{a\in \mathcal{A}}t_{a}(x_{a})x_{a}^{c}
\\
\label{eqn: total flow mixed time}
 \textbf{x}=\textbf{x}^{c}+\textbf{x}^{nc}   
\\
\label{eqn: total flow CAV on each link}
\textbf{x}^{c}=\textbf{A}\textbf{P}_{c}^{T}\textbf{g}^{c}
\\
\label{eqn: travel time func mixed}
t_{a}(x_{a})=t_{a}^{0}\sum _{i=1}^{n}\beta_{i}(\frac{x_{a}}{m_{a}})^{(i-1)}
\\
\label{eqn: probability constraint mixed CAV time}
\begin{array}
[c]{lr}
\sum_{r\in \mathcal{R}_{i}}p_{ir}^{c}=1;& \forall i\in [\![\mathcal{W}]\!]
\end{array}
\\
\begin{array}
[c]{lr}
 p_{ir}^{c}\in [0,1];&;\forall i\in [\![\mathcal{W}]\!],\forall r\in \mathcal{R}_{i}
\end{array}
\end{gather}
\end{subequations}
Constraint \eqref{eqn: total flow mixed time} states that the total flow in the network is the summation of the CAV flow ($\textbf{x}^c$) and  non-CAV flow ($\textbf{x}^{nc}$). Notice that we are minimizing the travel time for the CAV share of traffic. As a result, in  \eqref{eqn: mixed CAV time  objective} the traveling time of each link which is a function of both CAV flow and non-CAV flow ($t_a(x_a)$), is multiplied by the flow of CAVs only ($x_{a}^{c}$). The inputs to the optimization problem are the link-route incidence matrix $\textbf{A}$, O-D demand vector $\textbf{g}$, and non-CAV flow $\textbf{x}^{nc}$.

By solving Problem \ref{eqn: mixed time main problem}, we find optimal flows of CAVs over each O-D pair (route-probability matrix $\textbf{P}_c$). In other words, when a CAV enters the network at an origin \textit{O} given its destination \textit{D}, the algorithm gives it the desired socially optimal route to follow in terms of a sequence of links. 

As stated in Section 2.4 of  \cite{patriksson_traffic_2015} the system-centric problem can be reformulated as a user-centric problem by slightly changing the travel cost function. Therefore, the results on the existence and uniqueness of the solution for the user-centric problem (Section \ref{subsec: wardprop})  extend to the system-centric case. As a requirement for such a result we need  positive and strictly increasing travel time functions on $\textbf{P}$ which is achieved by having increasing polynomial functions.

\section{System-Centric Eco-Routing in the Presence of Mixed Traffic}\label{sec: mixed traffic energy}
In this section we solve the eco-routing problem for a fleet of CAVs in the presence of mixed traffic.  Eco-routing refers to the procedure of finding the optimal route for a vehicle to travel between two points which utilizes the least amount of energy costs. This problem shares similar properties to  \eqref{eqn: mixed time main problem}, with the difference that we minimize energy instead of time. As a result, we need an energy model to calculate energy consumption on each link. 

In this section, we first review an energy model for conventional vehicles and then formulate the system-centric eco-routing problem for CAVs.


\subsection{Energy Consumption Modeling} \label{sec: Empirical energy model}
\begin{figure}[h]
\centering
 \includegraphics[width=0.8\linewidth]{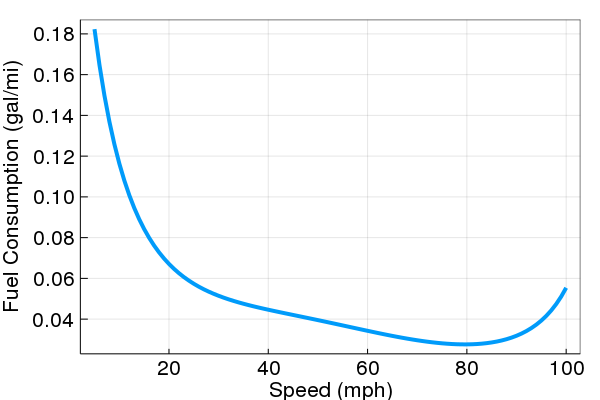}
  \caption{Average fuel consumption of conventional vehicles using Boriboonsomsin model (\(R_a=0\))}  
  \label{fig: fuel consumption CV}
\end{figure}

 Energy consumption of vehicles depends on many different factors including velocity and acceleration \cite{kamal_ecological_2011} of the vehicle, as well as the power-train's architecture.  Since in eco-routing we are making high-level decisions that can affect the energy consumption, a low-fidelity  model can be sufficient for our needs. Moreover, when solving the eco-routing problem, we are dealing with a large number of decision variables. Having a model with a simple mathematical function would allow us to speed up the calculation for practical purposes. Hence, we are looking for an energy model which can estimate the energy consumption as a function of the average speed of a vehicle. 
 We adopt the empirical energy model for conventional vehicles proposed by Boriboonsomsin et al. \cite{boriboonsomsin_eco-routing_2012}. This model is a polynomial function of of the average speeds of links. According to this empirical model (which is calibrated for an Audi A8), the average fuel consumption in grams per mile for every link $a\in \mathcal{A}$ can be calculated as follows:
\begin{equation}
\label{eqn: energy model CV}
    ln(e_{a})=\sum _{i=0}^{4}\theta_{i}(v_{a})^{i}+\theta_{5}R_{a}
\end{equation}
in which $e_{a}$ is the average energy consumption on link $a$ in $g/mi$, $v_{a}$ is the average speed of the link in \textit{mph}, $R_{a}$ is the road grade (in percentage), and \(\boldsymbol{\theta}=(\theta_{i}, i=0,1,...,5) \) is the vector of coefficients for calculating the energy cost. Typical values of $\boldsymbol{\theta}$ are given in Table \ref{tab: conversion factors energy CV}. Average fuel consumption per mile using  \eqref{eqn: energy model CV} and $\theta$ values in Table \ref{tab: conversion factors energy CV} is shown in Fig. \ref{fig: fuel consumption CV}. 
\begin{table}[h]
\caption{Energy cost coefficients \cite{boriboonsomsin_eco-routing_2012}%
}%
\label{tab: conversion factors energy CV}
\centering
\resizebox{1\columnwidth}{!}{\renewcommand{\arraystretch}{0.9}
\begin{tabular}{cccccc}
\hline
$\theta_{0}$ & $\theta_{1}$ & $\theta_{2}$ & $\theta_{3}$ & $\theta_{4}$ & $\theta_{5}$ \\ \hline
6.80         & -1.4e-1      & 3.92e-3      & -5.20e-5     & 2.57e-7      & 1.37e-1      \\ \hline
\end{tabular}
}\end{table}
In Appendix \ref{sec: CD/CS energy model} we also review a charge depleting(CD)/charge sustaining(CS) energy model \cite{karabasoglu_influence_2013} which can be used for PHEVs, HEVs, and EVs.
\subsection{Eco-routing Problem Formulation for Conventional Vehicles} \label{sec: eco-rputing conventional}
In order to formulate the eco-routing problem for conventional vehicles, we use energy model \eqref{eqn: energy model CV}. This problem is almost the same as \eqref{eqn: mixed time main problem}, with the only difference that $t_a(x_a)$ should be replaced with $e_a(x_a)$, which is the average fuel consumption per mile for traveling link $a\in \mathcal{A}$ as seen in \eqref{eqn: energy model CV}. Considering this, we rewrite the eco-routing problem of CAVs for conventional vehicles as follows:
 \begin{subequations}\label{eqn: mixed CAV main problem energy CV}
 \begin{gather}
\label{eqn: mixed CAV energy CV  objective}
\min_{\textbf{P}_c}\sum_{a\in \mathcal{A}}c_{gas}l_{a}e_{a}(v_{a}(x_a))x_{a}^{c}
\\
\label{eqn: total flow mixed energy CV}
 \textbf{x}=\textbf{x}^{c}+\textbf{x}^{nc}   
\\
\label{eqn: total flow CAV on each link energy CV}
\textbf{x}^{c}=\textbf{A}\textbf{P}_{c}^{T}\textbf{g}^{c}
\\
\label{eqn: travel time func mixed energy CV}
t_{a}(x_{a})=t_{a}^{0}\sum _{i=0}^{n}\beta_{i}(\frac{x_{a}}{m_{a}})^{i}
\\
\label{eqn: average speed energy CV}
v_{a}(x_{a})=\frac{l_{a}}{t_{a}(x_{a})}
\\
\label{eqn: energy model CV mixed CAV}
    ln(e_{a})=\sum _{i=0}^{4}\theta_{i}(v_{a})^{i}+\theta_{5}R_{a}
\\
\label{eqn: probability constraint mixed CAV time}
\begin{array}
[c]{lr}
\sum_{r\in \mathcal{R}_{i}}p_{ir}^{c}=1;& \forall i\in [\![\mathcal{W}]\!]
\end{array}
\\
\begin{array}
[c]{lr}
 p_{ir}^{c}\in [0,1];&\forall i\in [\![\mathcal{W}]\!],\forall r\in \mathcal{R}_{i}
\end{array}
 \end{gather}
 \end{subequations}
where $c_{gas}$ is the cost of gas (\$/gal), and $l_a$ is the length of link  \(\ a\in \mathcal{A}\). Moreover, $e_a$ is the average energy consumption per link's length \(\forall a\in \mathcal{A}\), and \(\boldsymbol{\theta}=(\theta_{i}, i=0,1,...,5) \) is the energy cost coefficient (Table \ref{tab: conversion factors energy CV}). 

The eco-routing problem formulation for PHEVs is formulated in Appendix \ref{sec: eco-route PHEV}.

\section{Non-CAV Flow Modeling} \label{sec: non-CAV flow modeling}
One of our assumptions is that the non-CAV flow is an input to our models, and can be inferred from actual traffic data. 
However, since we currently do not have CAVs in cities, we model the non-CAV flow by considering how non-CAVs react to the optimal decisions made by CAVs. To achieve this task, we assume  non-CAVs act selfishly by minimizing their travel time. This modeling framework has been extensively studied and is often referred to as the \emph{Traffic Assignment Problem} (TAP)\cite{patriksson_traffic_2015}. As a result, we propose an iterative method for finding  non-CAV flow considering the routing decisions of CAVs. The basis of this methodology is that whenever CAVs change their routing decisions, non-CAVs adjust theirs and vice versa. This process is well-known in game theory and is referred to as a Stackelberg game \cite{roughgarden_stackelberg_2004}.
For this particular problem, we consider an iterative procedure to find an equilibrium for mixed traffic flow of CAVs and non-CAVs. 

In order to obtain the non-CAV flow for a given CAV penetration rate $\gamma$, we first consider only  non-CAVs in the network and the O-D demand of non-CAVs is given by:
\begin{equation}
\textbf{g}^{nc}=(1-\gamma)\textbf{g}    
\end{equation}
Even though we choose a uniform demand distribution for non-CAVs between O-D pairs, without loss of generality, we can use any other given demand for both CAVs and non-CAVs. 

\begin{figure}[h]
\centering
\begin{subfigure}[b]{0.24\textwidth}
  \centering
  \includegraphics[width=1\columnwidth]{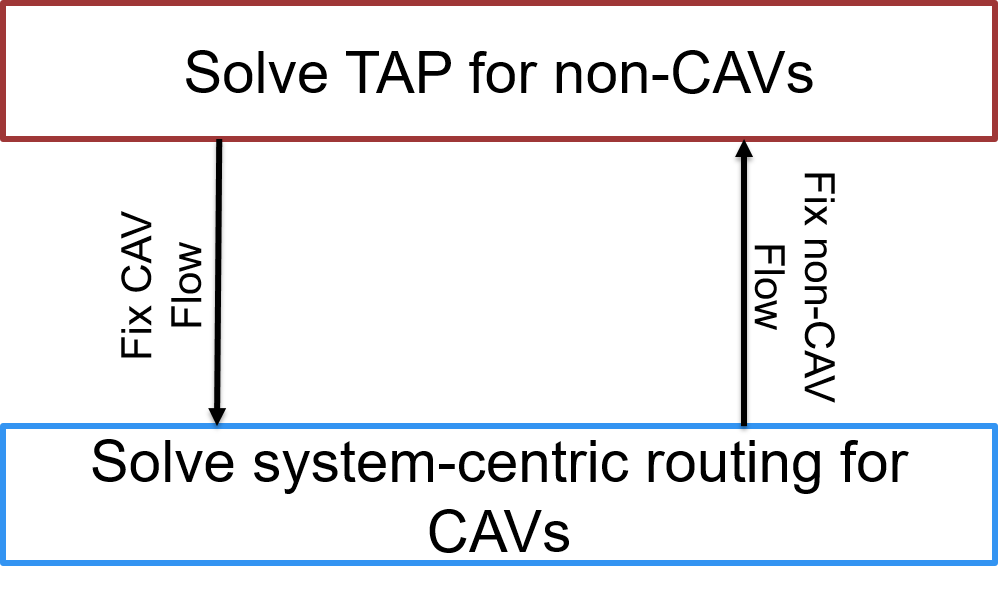}
  \caption{}
\label{fig: time-routing mixed CAV itteration}
\end{subfigure}%
\begin{subfigure}[b]{0.24\textwidth}
  \centering
  \includegraphics[width=1\columnwidth]{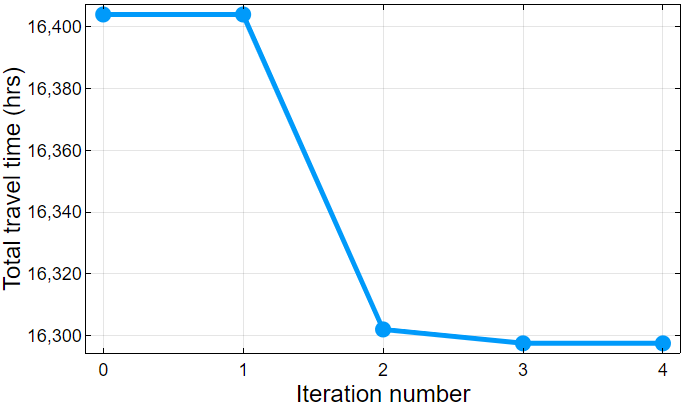}
  \caption{}
\label{fig: convergence plot for itterations-time}
\end{subfigure}
\caption{ (a) Procedure for solving the system-centric routing problem; (b) Convergence plot for iterating through TAP and social problem }
\label{fig: non-CAV flow iter process}
\end{figure}

Considering a non-CAV demand $\textbf{g}^{nc}$, we solve the selfish (user-centric) routing problem which minimizes their travel time. In this respect, we use the \emph{Method of Successive Averages} (MSA)\cite{sheffi_equilibrium_1985}.
After finding $\textbf{x}^{nc}$ using the  MSA, we solve the  time optimal (\ref{eqn: mixed time main problem}) or energy optimal (\ref{eqn: mixed CAV main problem energy CV}) routing problem for the CAV portion of traffic considering its demand to be: 
\begin{equation}
\textbf{g}^{c}=\gamma\textbf{g}
\end{equation}
Since non-CAVs were unaware of CAVs in the system while solving the TAP, we re-solve the problem considering CAV flow on each link. Hence, we re-iterate by considering the CAV solution $\textbf{x}^{c}$. Furthermore, the TAP is solved again for non-CAVs. Re-iteration of this process continuous until convergence (Figs. \ref{fig: non-CAV flow iter process}).
\subsection{Traffic Assignment Problem (TAP) and Wardrop equilibrium} \label{subsec: wardprop}

The objective of the Traffic Assignment Problem is to find link flows in a transportation network given the O-D demands and cost functions. A standard solution to this problem is to find  travel flows that minimize their travel times. Such a solution individually optimizes every vehicle's travel time based on network conditions.  This leads to a \emph{Nash Equilibrium} that in transportation networks is known as the \emph{Wardrop Equilibrium} \cite{beckmann_studies_1955}. The resulting flows $\mathbf{x^{*}}$ (equilibrium flows) require that for every O-D pair $\mathbf{w}$, and any route $r$ connecting $(w_{s},w_{t})$, the associated travel time is not greater than the traveling time from any other route. Formally
\begin{equation}
	t_a(x_a^*) \leq t_a(x_{a'}^*) \ \ \ \forall a,a' \in \mathcal{A}%
\end{equation}
or equivalently
\begin{equation}
t_a(x_r^*) \leq t_r(x_{r'}^*) \ \ \ \forall r,r' \in \mathcal{R}_i, \ \ \ \forall i \in[\![\mathcal{W}]\!] %
\end{equation}
To obtain such flows, we can solve
\begin{equation}\label{eqn: TAP original}
\min_{\mathbf{x} \in \mathcal{F}} \ \  \Phi(\mathbf{x}) = \sum\limits_{a \in \mathcal{A}}{\int\limits_{0}^{x_a}{t_a(s)ds}}
\end{equation}
where $\mathcal{F}$ is the set of feasible flow vectors defined by
\begin{equation}
\mathcal{F} = \Big\{  \mathbf{x}:\exists {\mathbf{x}^{\mathbf{w}}} \in \mathbb{R}_ +
		^{\left| \mathcal{A} \right|} ~\text{s.t.}~\mathbf{x} =
		\sum\limits_{\mathbf{w} \in \mathcal{W}} {{\mathbf{x}^{\mathbf{w}}}},\, \mathbf{N}{\mathbf{x}^{\mathbf{w}}} = {\mathbf{d}^{\mathbf{w}}},\,\forall \mathbf{w} \in \mathcal{W} \Big\}, 
\end{equation}
and where \(\mathbf{x}^{\mathbf{w}}\) is the flow vector attributed to O-D pair \(\mathbf{w}\). Recall that \(t_{a}(\cdot)\) in \eqref{eqn: travel time func mixed} is continuous. Since \(\mathcal{F}\) is a compact set,  the Weierstrass Theorem \cite{beckmann_studies_1955} implies that there exists a solution to this minimization problem. Moreover, since cost functions are non-decreasing (by assumption), then \(\Phi(\cdot)\) is convex and therefore a unique solution exists \cite{beckmann_studies_1955}.

Now, let us write the TAP in terms of non-CAV flows and take into account the presence of the CAV flow in the network. 
\begin{subequations}\label{eqn: TAP}
 \begin{gather}
\label{eqn: TAP objective}
\min_{\mathbf{x}^{nc}}  \sum\limits_{a \in \mathcal{A}}{\int\limits_{x^c_a}^{x^c_a + x_a^{nc}}{t_a(s)ds}}  
\\
s.t \ \  \mathbf{x}^{nc} = \sum\limits_{\mathbf{w} \in \mathcal{W}} {{\mathbf{x}^{nc, w}}} 
\\
\hspace{2,7cm} \mathbf{N}{\mathbf{x}^{nc, w}} = {\mathbf{d}^{nc, w}}, \ \  \forall \mathbf{w} \in \mathcal{W} \\
\mathbf{x}^{nc, w} \geq \textbf{0} \\ \notag
\end{gather}
\end{subequations}

\section{Numerical Results} \label{sec: numerical results}
In order to validate the proposed routing algorithms  we perform two case studies. First we analyze the widely explored Braess network (Fig. \ref{fig: Braess's Network}), and study the effect of the CAV penetration rate on the total time savings and energy savings in this network. As an alternative benchmark, we applied the algorithms to a sub-network of the Eastern Massachusetts interstate highways (Fig. \ref{fig:EMA-small-subnet}). For finding the energy optimal routes, we assume the road grade is zero (\(R_{a}=0\) in  \eqref{eqn: energy model CV}), and we assume the cost of gas is  2.75 \$/gal. We solve the NLP problems using IPOPT \cite{wachter_implementation_2006} in Julia \cite{bezanson_julia:_2017}.

For the eco-routing case, we only show the results for solving  \eqref{eqn: mixed CAV main problem energy CV}. In other words, we only solve the eco-routing problem for conventional vehicles using the energy model discussed in Section \ref{sec: Empirical energy model}. 
As mentioned before, eco-routing results are extremely sensitive to the energy model. Given a more accurate energy model which is convex, smooth and differentiable we may get different results. Hence, the eco-routing results shown in this paper should only be considered as preliminary results which show the potential of saving energy using centralized routing of CAVs.

\subsection{Braess Network Example}
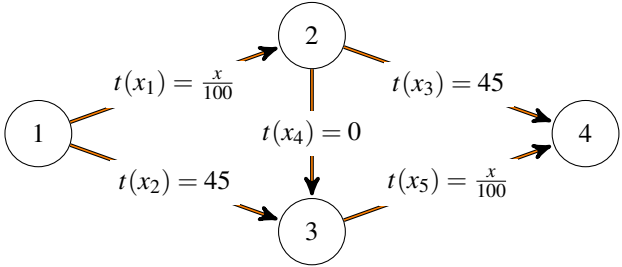
\begin{figure}[h]
\centering
\begin{tikzpicture}[>=stealth',shorten >=1pt,node distance=3.3cm,on grid,initial/.style    ={}]
  \node[state]          (1)                        {$1$};
  \node[state]          (2) [above right =1.3cm and 3.6cm of 1]    {$2$};
  \node[state]          (3) [below right =1.3cm and 3.6cm of 1]    {$3$};
  \node[state]          (4) [below right =1.3cm and 3.6cm of 2]    {$4$};
\tikzset{mystyle/.style={->,double=orange}} 
\tikzset{every node/.style={fill=white}} 
\path (2)     edge [mystyle]    node   {$t(x_{3})=45$} (4)
      (1)     edge [mystyle]    node   {$t(x_{1})=\frac{x}{100}$} (2);
\tikzset{mystyle/.style={->,double=orange}}   
\path (1)     edge [mystyle]   node   {$t(x_{2})=45$} (3)
      (3)     edge [mystyle]   node   {$t(x_{5})=\frac{x}{100}$} (4) 
      (2)     edge [mystyle]   node   {$t(x_{4})=0$} (3) ;

\end{tikzpicture}
\caption{Simple 5 link directed graph (Braess's network)}
\label{fig: Braess's Network}
\end{figure}

\begin{figure}
        \begin{subfigure}{.49\columnwidth}
          \centering
          \includegraphics[width=1\columnwidth]{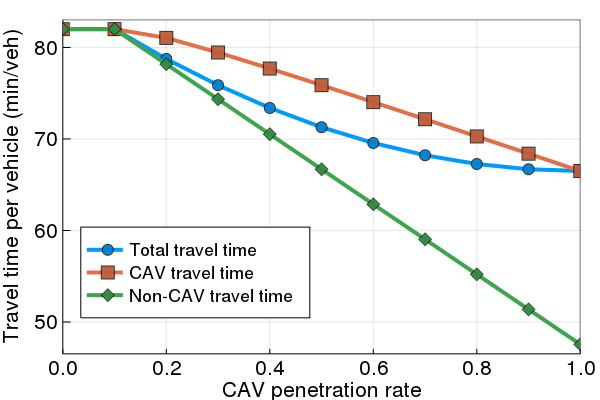}
          \caption{}
          \label{fig: Braess graph penetration rate time opt- Time}
        \end{subfigure} 
        \begin{subfigure}{.49\columnwidth}
          \centering
          \includegraphics[width=1\columnwidth]{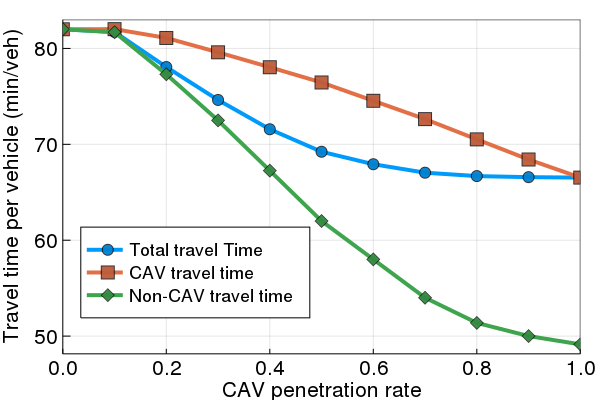}
          \caption{}
          \label{fig: Braess graph penetration rate energy opt- Time}
        \end{subfigure}
        \begin{subfigure}{.49\columnwidth}
          \centering
          \includegraphics[width=1\columnwidth]{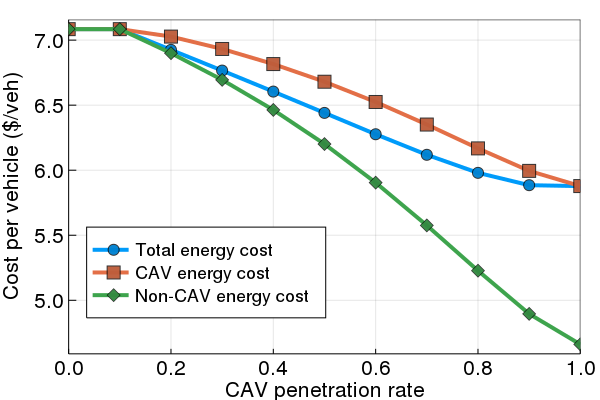}
          \caption{}
          \label{fig: Braess graph penetration rate time opt-energy}
        \end{subfigure} 
        \begin{subfigure}{.49\columnwidth}
          \centering
          \includegraphics[width=1\columnwidth]{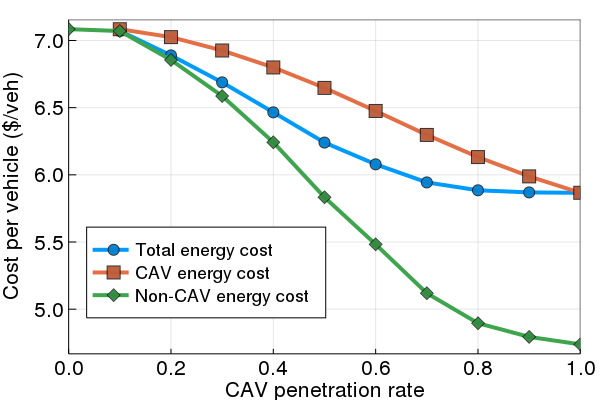}
          \caption{}
          \label{fig: Braess graph penetration rate energy opt-energy}
        \end{subfigure}
        \begin{subfigure}{.49\columnwidth}
          \centering
          \includegraphics[width=1\columnwidth]{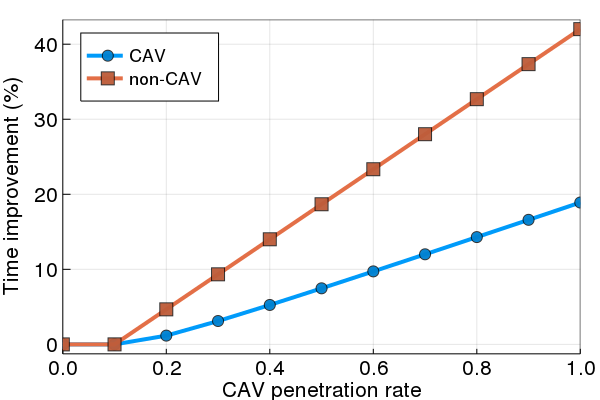}
          \caption{}
          \label{fig: Braess graph penetration rate time opt-improvements}
        \end{subfigure} 
        \begin{subfigure}{.49\columnwidth}
          \centering
          \includegraphics[width=1\columnwidth]{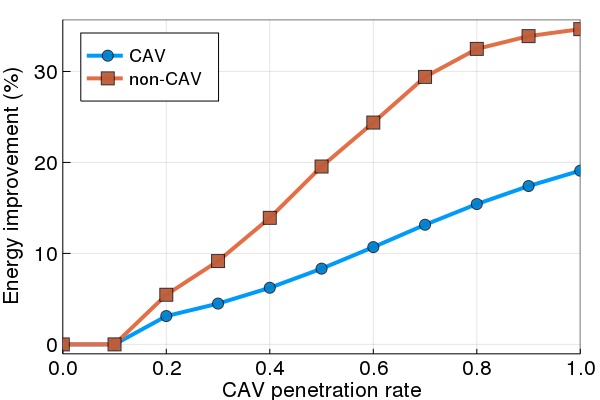}
          \caption{}
          \label{fig: Braess graph penetration rate energy opt-improbements}
        \end{subfigure}
        \caption{Braess network routing results under different penetration rates by solving system-centric time-optimal problem (a), (c), (e) and, system-centric energy-optimal problem (b), (d), (f).}
        \label{fig: Braess energy optimal plots} 
    \end{figure}

To demonstrate  how  optimal routing of CAVs under different  penetration rates can affect both CAVs and non-CAVs, we first apply algorithm \ref{eqn: mixed time main problem} to the well-known Braess network (Fig. \ref{fig: Braess's Network}).  Note that in this case instead of using the BPR function \eqref{eqn: travel time func mixed}, we use the travel time functions shown on each link of the Braess network in Fig. \ref{fig: Braess's Network}.  We consider a demand of 4000 veh/hr travels from node 1 to node 4, the lengths of links 1,2,3 and 5 equal to 30.5 miles and the length of link 4 equal to zero.  First we solve the time-optimal routing of CAVs under different penetration rates. Using the obtained flows, and energy model \eqref{eqn: energy model CV} we calculate energy costs for traveling through the network (all cars are assumed to be conventional vehicles). Time-optimal results are shown in Figs. \ref{fig: Braess graph penetration rate time opt- Time}, in which we compare traveling time of CAVs with non-CAVs under different penetration rates. The energy cost for traveling through the optimal routes are also shown in Fig. \ref{fig: Braess graph penetration rate time opt-energy}. In addition, we compared the traveling time of CAVs and non-CAVs under different penetration rates using the case of 0\% CAV as a baseline and reported the time savings in Fig. \ref{fig: Braess graph penetration rate time opt-improvements}. As shown in Fig. \ref{fig: Braess energy optimal plots}, introducing CAVs into the system not only improves the time saving of CAVs, but also helps non-CAVs to save time. This is because smart routing decisions of CAVs reduce the traffic intensity in the highly congested roads, which consequently helps non-CAVs to travel faster. Note that the baseline is the 0 \% CAV case (uncontrolled traffic), in which all vehicles act selfishly. As we inject CAVs into the system, we see that travel time (as well as energy cost) per vehicle of CAVs starts decreasing compared with the uncontrolled traffic. Moreover, the traveling time of commuting through the fastest route decreases as we inject more CAVs to the system. Typically, we expect a trade off between time saving and energy saving in routing problems \cite{sun_save_2016, houshmand_eco-routing_2018}. However, in Fig. \ref{fig: Braess energy optimal plots} we see that time and energy follow the same trend. In other words, time savings result in energy savings. The reason for this behavior is the energy model used in the eco-routing problem \ref{eqn: energy model CV}. As we can see in Fig. \ref{fig: fuel consumption CV}, the higher the speed, the better fuel efficiency. Hence, for conventional vehicles, and based on this energy model \cite{boriboonsomsin_eco-routing_2012}, we get similar results for energy and time.

It is interesting to see that when a small percentage of CAVs are in the system, there is no improvement for anyone. This happens because CAVs are optimizing over their own small fraction of overall traffic and this fraction is not sufficient to change the network conditions. However, as we increase the penetration rate, CAVs create a more balanced flow distribution in the network from which both CAVs and non-CAVs can benefit. In the Braess network example, it can be seen that if all the cars in the system are replaced with CAVs, we can save 18.9\% in terms of travel time. This value is often referred to as the Price of Anarchy (PoA) \cite{zhang_price_2018}.

In addition to the time-optimal case, we also solve the eco-routing (energy-optimal) problem for CAVs using the Braess network. As mentioned earlier, there are many different  models to calculate energy costs of vehicles which depend on the vehicle type. In Section \ref{sec: mixed traffic time} we formulated the eco-routing routing problem using an empirical energy model for conventional vehicles in  \eqref{eqn: mixed CAV main problem energy CV}. In Appendix \ref{sec: eco-route PHEV}, we formulate the system-centric eco-routing problem for PHEVs as shown in \eqref{eqn: mixed CAV energy PHEV}. However, due to the non-convexity of \eqref{eqn: mixed CAV energy PHEV} and issues with local optima, we are only showing the energy-optimal results by solving \eqref{eqn: mixed CAV main problem energy CV}. 
As we can see in Figs. \ref{fig: Braess graph penetration rate energy opt- Time}, \ref{fig: Braess graph penetration rate energy opt-energy}, and \ref{fig: Braess graph penetration rate energy opt-improbements}, energy-optimal results follow the same trend as time-optimal result. In other words, centralized eco-routing of CAVs can benefit both CAVs and non-CAVs. The maximum energy savings happens at the 100\% CAV penetration rate (19.1\%).

\subsection{EMA Interstate Highway Network }
\begin{figure}[h]
\centering
\begin{subfigure}[b]{0.2\textwidth}
  \centering
  \includegraphics[width=1\columnwidth]{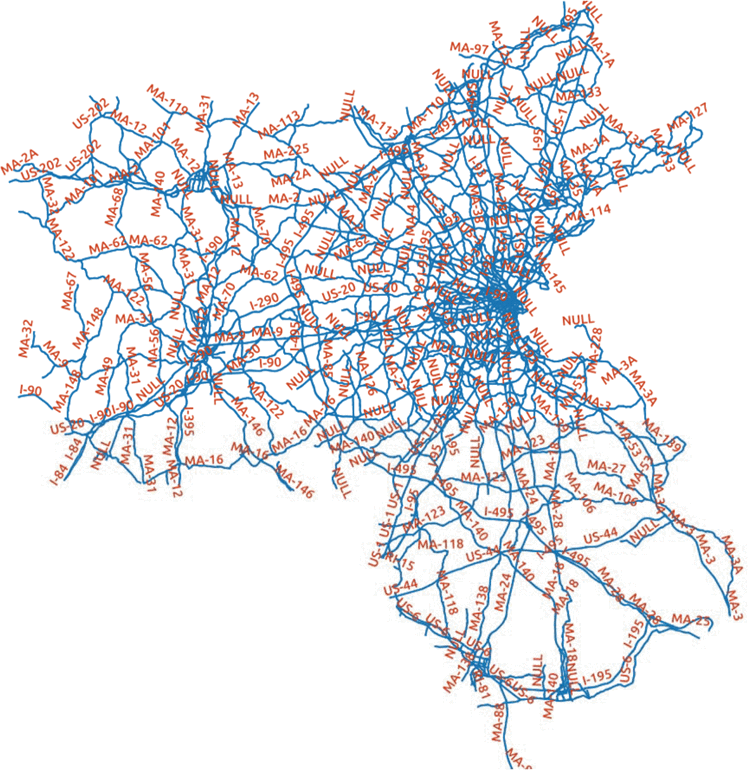}
  \caption{}
  \label{fig:EMA-All-TMC}
\end{subfigure}%
\begin{subfigure}[b]{0.2\textwidth}
  \centering
  \includegraphics[width=1\columnwidth]{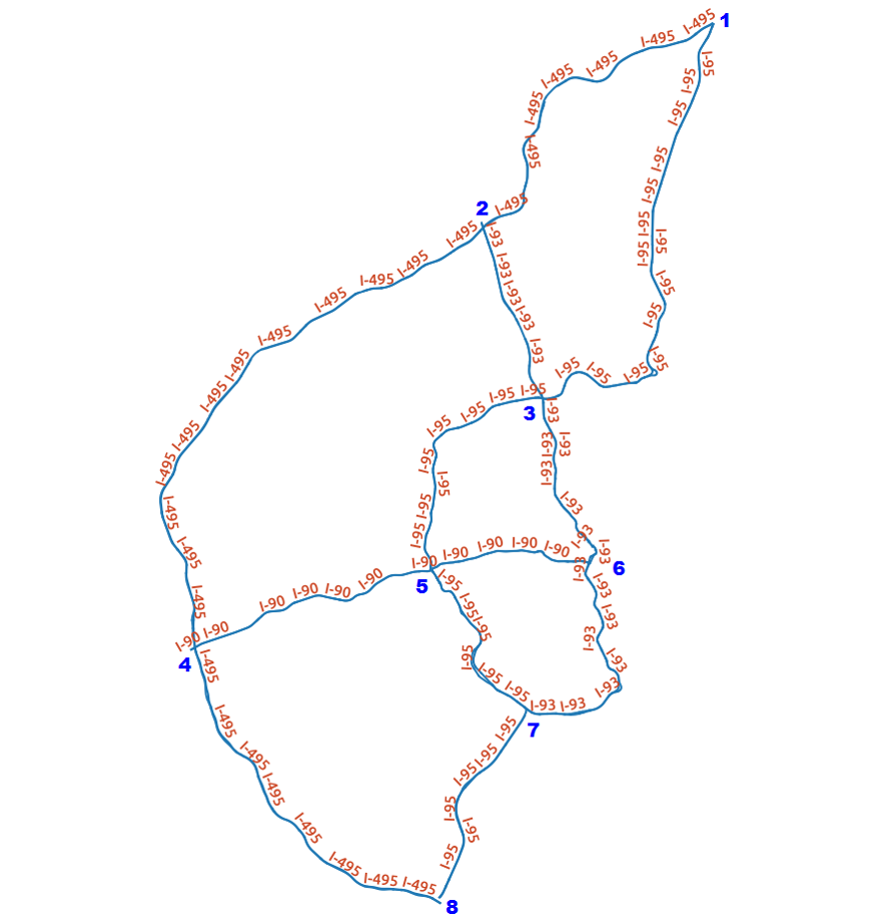}
  \caption{}
  \label{fig:EMA-small-subnet}
\end{subfigure}
\caption{ (a) All available road segments in the road map of Eastern Massachusetts \cite{zhang_price_2018}  ; (b) Interstate highway sub-network of eastern Massachusetts}
\label{C}
\end{figure}

\begin{figure}[h]
\centering
\begin{subfigure}{.49\columnwidth}
  \centering
  \includegraphics[width=1\columnwidth]{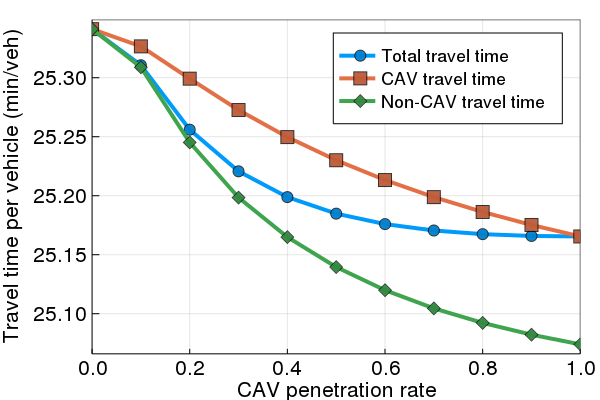}
  \caption{}
  \label{fig: EMA graph penetration rate time opt- Time}
\end{subfigure}
\begin{subfigure}{.49\columnwidth}
  \centering
  \includegraphics[width=1\columnwidth]{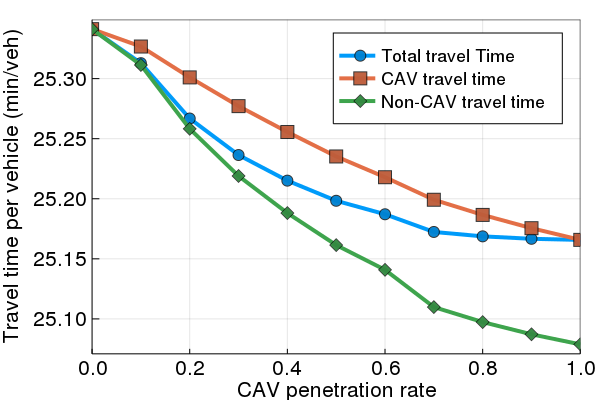}
  \caption{}
  \label{fig: EMA graph penetration rate energy opt- Time}
\end{subfigure}
\begin{subfigure}{.49\columnwidth}
  \centering
  \includegraphics[width=1\columnwidth]{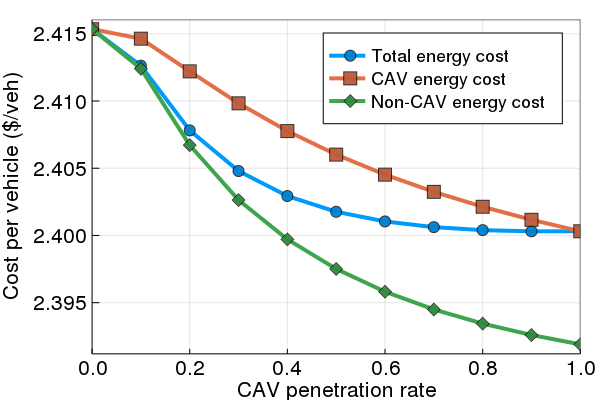}
  \caption{}
  \label{fig: EMA graph penetration rate time opt-energy}
\end{subfigure}%
\begin{subfigure}{.49\columnwidth}
  \centering
  \includegraphics[width=1\columnwidth]{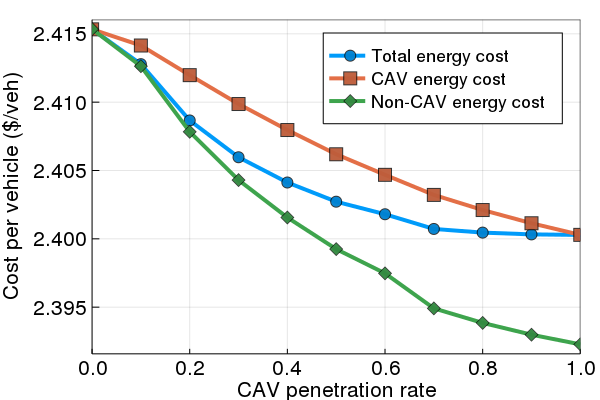}
  \caption{}
  \label{fig: EMA graph penetration rate energy opt-energy}
\end{subfigure}
\begin{subfigure}{.49\columnwidth}
  \centering
  \includegraphics[width=1\columnwidth]{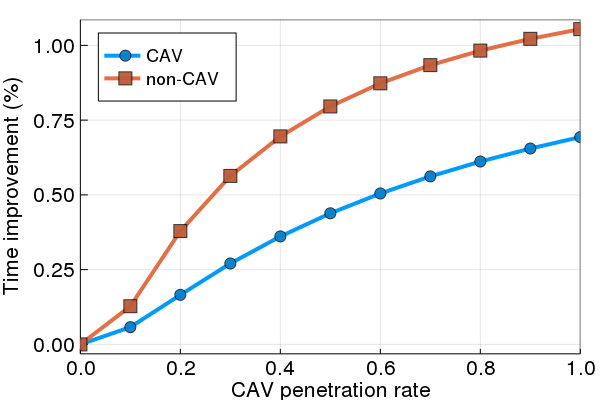}
  \caption{}
  \label{fig: EMA graph penetration rate time opt-improvements}
\end{subfigure}
\begin{subfigure}{.49\columnwidth}
  \centering
  \includegraphics[width=1\columnwidth]{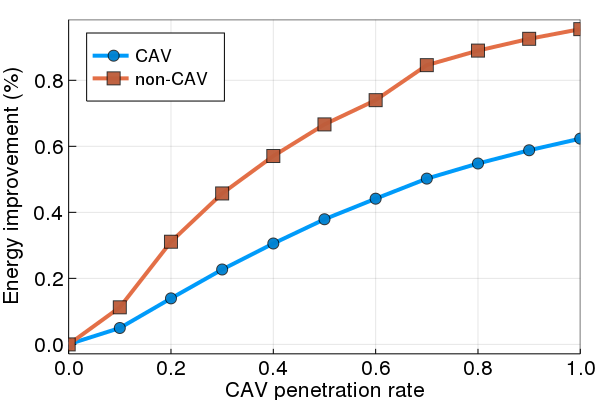}
  \caption{}
  \label{fig: EMA graph penetration rate energy opt-improbements}
\end{subfigure}
\caption{EMA network routing results under different penetration rates by solving system-centric  time-optimal  problem  (a), (c), (e) and system-centric  energy-optimal  problem  (b), (d), (f).}
\label{fig: EMA time optimal plots}
\end{figure}

\begin{figure}[h]
\centering
\begin{subfigure}[b]{0.24\textwidth}
  \centering
  \includegraphics[width=1\columnwidth]{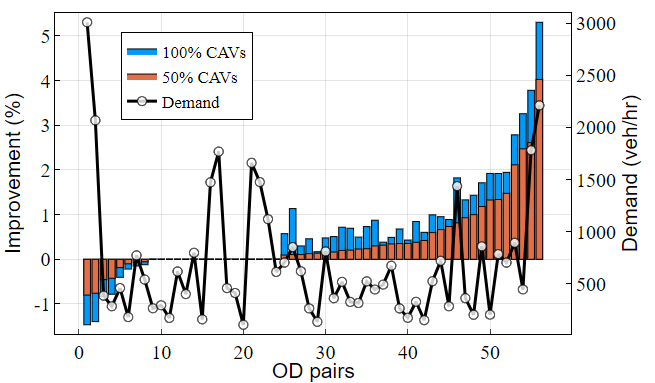}
  \caption{}
\label{fig: imrovement time 50 and 100}
\end{subfigure}%
\begin{subfigure}[b]{0.24\textwidth}
  \centering
  \includegraphics[width=1\columnwidth]{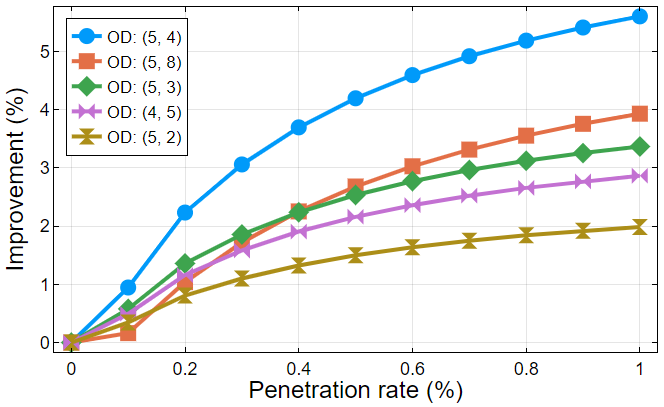}
  \caption{}
\label{fig:improvement time top 5}
\end{subfigure}
\caption{ (a) Time saving improvement of different O-D pairs for 50\% and 100\% CAV time optimal routing (EMA network)  ; (b) Travel time improvement of different O-D pairs as a function of CAV penetration rate (EMA network)}
 \label{fig:EMA Interstate improvements}
\end{figure}

In order to obtain more realistic results, we perform a data-driven case study using the actual
traffic data from the Eastern Massachusetts (EMA) road network. These data were collected by \textit{INRIX} and provided to us by the Boston Region
Metropolitan Planning Organization. The sub-network including
the interstate highways of EMA (Fig. \ref{fig:EMA-small-subnet}) is chosen for
the case study. For this network, we use the O-D demand which has been estimated using an inverse optimization framework in \cite{zhang_price_2018}. In order to solve the problem we consider 56 O-D pairs, and allow up to 3 routes between each origin and destination (top 3 shortest routes). We then solve  \eqref{eqn: mixed time main problem} and \eqref{eqn: mixed CAV main problem energy CV} in order to find the time optimal and energy optimal paths for CAVs respectively. Time optimal results are shown in Fig. \ref{fig: EMA time optimal plots}, and the energy optimal results are shown in Fig. \ref{fig: Braess energy optimal plots}. The results follow the same behavior as the results of the Braess example. We again see that as the CAV penetration rate increases, both CAVs and non-CAVs benefit from optimal routing decisions of non-CAVs. In Fig. \ref{fig: imrovement time 50 and 100}, we show the time improvement of different O-D pairs with their corresponding O-D demand for 100\% and 50\% CAV penetration rates. It is interesting to see in Fig. \ref{fig: imrovement time 50 and 100} that two OD pairs with relatively high demand are being affected by -0.8\% and -1.3\% for 50\% and 100\% $\gamma$'s respectively. However, we see improvements on most of the OD pairs. In this manner, we are able to identify which OD pairs are getting worse and which ones are improving. This gives the opportunity to better understand the dynamics of the  network.

\section{Conclusions and Future Work} \label{sec: conclusions}
In this paper we proposed system-centric optimal routing algorithms for a fleet of CAVs in the presence of mixed traffic. We consider two objectives for routing: (1) minimizing travel time (2) minimizing energy consumption cost. Moreover, in order to model the routing behavior of regular vehicles, we assume that they make selfish decisions by minimizing their own travel time. Then, by iteratively solving the TAP, and finding optimal routes for CAVs we estimate the non-CAV flow in the network. Historical traffic data and a simple illustrative example were used to validate the models. The results indicate that optimal routing of CAVs can not only benefit CAVs, but the smart routing decision of CAVs helps ease traffic congestion in the network which helps regular vehicles as well. Additionally, we empirically showed that even a small CAV penetration rate has significant impact on the overall traveling cost of the network.

So far, we assumed that all CAVs belong to the same fleet and we can route all of them. In ongoing work, we are considering multiple fleets of CAVs, in which each fleet is trying to minimize its own cost. Moreover, as  pointed out earlier, eco-routing of CAVs is highly dependent on the energy model. We are actively looking for other energy models which can predict the cost for different vehicle types (e.g., EVs, PHEVS) with reasonable levels of accuracy.

\bibliographystyle{IEEEtran}
\begin{tiny}
\bibliography{NEXTCAR}
\end{tiny}
\appendix 
\subsection{CD/CS Energy Model} \label{sec: CD/CS energy model}
Karabasoglu et al. \cite{karabasoglu_influence_2013} proposed an approximate energy model for different classes of vehicles including EVs, HEVs, PHEVs, and CVs. They considered two operational modes for vehicles: charge depleting (CD), and charge sustaining (CS). The CD mode refers to the case in which the main propulsion energy for driving the car comes from the battery pack (electricity). In addition, the CS mode occurs when the vehicle uses the internal combustion engine (gas) to drive the vehicle. PHEVs use both CD and CS modes since they have both engine and electric motor. HEVs only use the CS mode, and EVs only operate on the CD mode. Even though conventional vehicles do not have a battery pack and the CS mode does not apply to them, for the sake of consistency, we use the CS mode for conventional vehicles to refer to their gas operational mode. Karabasoglu et al. calculated the average $mi/gal$ and $mi/kWh$ that a vehicle can travel through different standard drive cycles (NYC, UDDS, HWFET, etc.) under CD or CS operational modes. They made these calculations using \textit{PSAT}, a commercially available software package simulating different power-train architectures using high-fidelity models. Calculated average energy consumption values under CD and CS modes are referred to as $\mu_{CD}$ and $\mu_{CS}$ respectively and are reported in Table \ref{tab: conversion factors}. Qiao et al. \cite{qiao_vehicle_2016} and Houshmand et al. \cite{houshmand_eco-routing_2018}  solved the eco-routing problem for PHEVs by benefiting from the values calculated in \cite{karabasoglu_influence_2013} (Table \ref{tab: conversion factors}).
\begin{table}[h]
\caption{Average Energy consumption values \cite{karabasoglu_influence_2013}%
}%
\label{tab: conversion factors}
\centering
\resizebox{1\columnwidth}{!}{\renewcommand{\arraystretch}{0.9}
\begin{tabular}{llllll}
\hline
Vehicle Type & Symbol     & Unit     & HWFET & UDDS & NYC  \\ \hline
PHEV         & $\mu_{CD}$ & $mi/kWh$ & 5.7   & 6.2  & 4.2  \\
             & $\mu_{CS}$ & $mi/gal$ & 58.6  & 69.4 & 45.7 \\ 
HEV          & $\mu_{CS}$ & $mi/gal$ & 59.7  & 69.5 & 48.0 \\ 
EV           & $\mu_{CD}$ & $mi/kWh$ & 5.2   & 4.8  & 3.1  \\ 
CV           & $\mu_{CS}$ & $mi/gal$ & 52.8  & 32.1 & 16.4 \\ \hline
\end{tabular}
}\end{table}
\subsection{Eco-routing Problem Formulation for PHEVs} \label{sec: eco-route PHEV}
In order to solve the eco-routing problem for CAVs, we follow the same formulation used in \cite{houshmand_eco-routing_2018}. The essence of this eco-routing model is that it categorizes each link based on its average speed into  3 different modes: heavy traffic, medium traffic, and low traffic links (note that we can have more than 3 modes as well). We then assign a drive cycle to each link based on its average speed (Table \ref{tab: Drive cycle assignment}). 
\begin{table}[h]
\caption{Drive cycle assignment of each link }%
\label{tab: Drive cycle assignment}
\centering
\resizebox{0.7\columnwidth}{!}{\renewcommand{\arraystretch}{0.9}
\begin{tabular}{ccc}
\hline
Traffic Mode   & \begin{tabular}[c]{@{}c@{}}Average speed \\ on the link (mph)\end{tabular} & Drive Cycle                              \\ \hline
Heavy Traffic    & [0,20]                                                                     & NYC         \\
Medium Traffic & [20,40]                                                                    & UDDS\\
Low Traffic   & $\geq$ 40                                                                    & HWFET                       \\ \hline
\end{tabular}
}\end{table}

We consider two sets of decision variables: the CAV route-probability matrix,  \(\textbf{P}_c=\{p_{ir}^c, i\in [\![\mathcal{W}]\!], r\in \mathcal{R}\}\) and the CD/CS switching strategy on each link, \(\textbf{Y}=\{y_{a,r}^{i}, i\in [\![\mathcal{W}]\!], r\in \mathcal{R}, a\in \mathcal{A} \}\). Here, $y_{a,r}^{i}$ represents the fraction of link length (\(l_a\))  over which we use the CD
mode. Thus, if we only use the CD mode over link \(a\in \mathcal{A}\), then \(y_{a,r}^{i}=1\). Note that in order to solve this problem for EVs, we should set $\textbf{Y}=1$; and for HEVs and CVs we let $\textbf{Y}=0$. Consequently we remove the dependency of the problem on $\textbf{Y}$ in these cases.
 
This problem shares similarities with \eqref{eqn: mixed time main problem}, with the only difference that here we are minimizing energy instead of time. Using the same notation introduced in sections \ref{sec: problem formulation 100 CAV time}, and \ref{sec: mixed traffic time}, for a given CAV penetration rate $\gamma$, we can formulate the problem as follows:

 \begin{subequations} \label{eqn: mixed CAV energy PHEV}
 \begin{gather}
\label{eqn: mixed CAV system-centric energy}
\begin{split}
\min_{\textbf{P}_{c},\textbf{Y}} & \sum_{i\in [\![\mathcal{W}]\!]}\sum_{r\in \mathcal{R}_{i}}\sum_{a\in \mathcal{A}}(c_{gas}\frac{l_{a}}{\mu_{cs}^{a}(v_a(x_a))}(1-y_{a,r}^{i}) \\ 
  & \hspace{3cm}+ c_{ele}\frac{l_{a}}{\mu_{CD}^{a}(v_a(x_a))}y_{a,r}^{i})x_{a}^{c}
\end{split}
\\
\label{eqn: battery constraint}
\begin{array}
[c]{lccc}%
s.t. & \sum_{a\in \mathcal{A}}\frac{\alpha_{a,r}^{i}y_{a,r}^{i}l_{a}}{\mu_{cd}^{a}}\leq E_{0}^{r,i}
& ;\forall r\in \mathcal{R}_{i}, & \forall i\in [\![\mathcal{W}]\!]
\end{array}
\\
\label{eqn: total flow mixed energy}
 \textbf{x}=\textbf{x}^{c}  +\textbf{x}^{nc}   
\\
\textbf{x}^{c}=\textbf{A}\textbf{P}_{c}^{T}\textbf{g}^{c}
\\
\label{eqn: average speed energy CD/CS}
v_{a}(x_{a})=\frac{l_{a}}{t_{a}(x_{a})}
\\
\begin{array}
[c]{lrr}
\textbf{P}_{c}=[p_{ir}^{c}] &, p_{ir}^{c}\in [0,1];&\forall i\in [\![\mathcal{W}]\!],\forall r\in \mathcal{R}_{i}
\end{array}
\\
\label{eqn: probability const. Social energy}
\begin{array}
[c]{lr}
\sum_{r\in \mathcal{R}_{i}}p_{ir}^{c}=1;& \forall i\in [\![\mathcal{W}]\!]
\end{array}
\\
\label{eqn: Y constraint mixed energy}
\begin{array}
[c]{lr}
y_{a,r}^{i}\in [0,1];&\forall a\in \mathcal{A}, \forall i\in [\![\mathcal{W}]\!],\forall r\in \mathcal{R}_{i}
\end{array}
\end{gather}
 \end{subequations}
 where $c_{gas}$ and $c_{ele}$ are the cost of gas (\$/gal), and electricity (\$/kWh) respectively.
 Constraint \eqref{eqn: battery constraint}, is the energy constraint stating that the total electrical energy used on each path and O-D pair should not be more than the available energy in the battery pack at the start of that path ($E_{0}^{r,i}$). $mu_{CS}$ and $mu_{CD}$ are functions of velocity (Table \ref{tab: conversion factors}, and \ref{tab: Drive cycle assignment}), and velocity is a function of flow \eqref{eqn: average speed energy CD/CS}. As mentioned in Sec. \ref{sec: mixed traffic time}, we assume that the  non-CAV traffic flow equilibrium is inferred from data, and is known  ($\textbf{x}_{nc}$).
By solving \eqref{eqn: mixed CAV energy PHEV}, when a vehicle enters the network at an origin \textit{O} given its destination \textit{D}, the algorithm gives it the desired socially optimal route to follow in terms of a sequence of links, and the optimal CD/CS switching strategy on each link.

This model finds the eco-route for PHEVs, as well as HEVs, and EVs by setting \(Y=0\) or $1$ respectively. We should note that the  energy model used in \eqref{eqn: mixed CAV energy PHEV} is a piece-wise constant function of velocity, and can cause non-convexity and differentiability issues in the problem solution. 
\end{document}